%%%%%%%%%%%%%%%%%%%%%%%%%%%%%%%%%%%%%%%%%%%%%%%%%%%%%%%%%%%%%%%
%%%%%%%%%%%%%%%%%%%%%%%%%%%%%%%%%%%%%%%%%%%%%%%%%%%%%%%%%%%%%%%

\magnification=1200
\overfullrule=0pt

\font\eightrm=cmr8

\textfont0=\tenrm
\scriptfont0=\sevenrm
\scriptscriptfont0=\fiverm

%%%%%%%%%%%%%%%%
%%  FONT CMM  %%
%%%%%%%%%%%%%%%%
%
% \font\sixi=cmmi6
% \font\eighti=cmmi8
% \font\ninei=cmmi9
% \font\twelvei=cmmi12
% \font\fourteeni=cmmi10 at 14pt
%
\font\teni=cmmi10
\font\seveni=cmmi7
\font\fivei=cmmi5
\textfont1=\teni
\scriptfont1=\seveni
\scriptscriptfont1=\fivei

%%%%%%%%%%%%%%%%%
%%  FONT CMSY  %%
%%%%%%%%%%%%%%%%%
%

%%%%%%%%%%%%%%%%%
%%  FONT CMTI  %%
%%%%%%%%%%%%%%%%%
%

\font\eightit=cmti8

%%%%%%%%%%%%%%%%%%%
%%  FONT CMBFTI  %%
%%%%%%%%%%%%%%%%%%%
%

%%%%%%%%%%%%%%%%%%%%%%%%%%%%%%%%%%%%%
%%  FONT CMCSC (Petites capitales) %%
%%%%%%%%%%%%%%%%%%%%%%%%%%%%%%%%%%%%%

\font\sc=cmcsc10

%%%%%%%%%%%%%%%%%
%%  FONT MSBM  %%
%%%%%%%%%%%%%%%%%

% \font\ninebb=msbm9  % msbm10 at 9pt
% \font\eightbb=msbm8 % msbm10 at 8pt
% \font\sixbb=msbm6   % msbm10 at 6pt
% \font\sixbboard=msbm7 at 6pt
%
\font\tenbb=msbm10
\font\sevenbb=msbm7
\font\fivebb=msbm5
\newfam\bbfam
\textfont\bbfam=\tenbb
\scriptfont\bbfam=\sevenbb
\scriptscriptfont\bbfam=\fivebb
\def\bb{\fam\bbfam\tenbb}
\let\oldbb=\bb
\def\bb #1{{\oldbb #1}}

%%%%%%%%%%%%%%%%%%%%%%%%%%%%%
%%  FONT EUFM  (Gothique)  %%
%%%%%%%%%%%%%%%%%%%%%%%%%%%%%
%

%
\font\tengoth=eufm10
\font\sevengoth=eufm7
\font\fivegoth=eufm5
\newfam\gothfam
\textfont\gothfam=\tengoth
\scriptfont\gothfam=\sevengoth
\scriptscriptfont\gothfam=\fivegoth

%

%%%%%%%%%%%%%%%%%
%%  FONT CMBX  %%
%%%%%%%%%%%%%%%%%
%              
% \font\eightbf=cmbx10 at 8pt

\font\fourteenbf=cmbx10 at 14pt

\font\eightbf=cmbx8
% \font\sevenbf=cmbx10 at 7pt
 
% \font\sixbf=cmbx5
% \font\fivebf=cmbx10 at 5pt

\font\tenbf=cmbx10
\font\sevenbf=cmbx7
\font\fivebf=cmbx5
\newfam\bffam
\textfont\bffam=\tenbf
\scriptfont\bffam=\sevenbf
\scriptscriptfont\bffam=\fivebf
\def\bf{\fam\bffam\tenbf}
%
%

%%%%%%%%%%%%%%%%%%%%%%%%%%%%%%%%%%%%%%%%%%%%%%%%%%%%%%%%%%%%

%%%%%%%%%%%%%%%%%%%%%%%%%%%%%%%%%%%%%%%%%%%%%%%%%%%%%%%%%%%%
%%  HAUT-DE-PAGE (EX. : AUTEUR COURANT ET TITRE COURANT)  %%
%%  BAS-DE-BAGE                                           %%
%%%%%%%%%%%%%%%%%%%%%%%%%%%%%%%%%%%%%%%%%%%%%%%%%%%%%%%%%%%%

\newif\ifpagetitre           \pagetitretrue
\newtoks\hautpagetitre       \hautpagetitre={\hfil}
\newtoks\baspagetitre        \baspagetitre={\hfil}

\newtoks\auteurcourant       \auteurcourant={\hfil}
\newtoks\titrecourant        \titrecourant={\hfil}

\newtoks\hautpagegauche
         \hautpagegauche={\hfil\the\auteurcourant\hfil}
\newtoks\hautpagedroite
         \hautpagedroite={\hfil\the\titrecourant\hfil}

\newtoks\baspagegauche
         \baspagegauche={\hfil\tenrm\folio\hfil}
\newtoks\baspagedroite
         \baspagedroite={\hfil\tenrm\folio\hfil}

\headline={\ifpagetitre\the\hautpagetitre
\else\ifodd\pageno\the\hautpagedroite
\else\the\hautpagegauche\fi\fi }

\footline={\ifpagetitre\the\baspagetitre
\global\pagetitrefalse
\else\ifodd\pageno\the\baspagedroite
\else\the\baspagegauche\fi\fi }

%%%%%%%%%%%%%%%%%%%%%%%%%%%%%%%%%%%%%%%%%%%%%%%%%%%%%%%%%%%%%%%

           %%%%%%%%%%%%%%%%%%%%%%%%%%%%%%%%%%%%%%%%%%%
           %%  UN SEMBLANT DES GUILLEMETS FRANCAIS  %%
           %%%%%%%%%%%%%%%%%%%%%%%%%%%%%%%%%%%%%%%%%%%

\def\og{\leavevmode\raise.3ex
     \hbox{$\scriptscriptstyle\langle\!\langle$~}}
\def\fg{\leavevmode\raise.3ex
     \hbox{~$\!\scriptscriptstyle\,\rangle\!\rangle$}}
     
%%%%%%%%%%%%%%%%%%%%%%%%%%%%%%%%%%%%%%%%%%%%%%%%%%%%%%%%%%%%%%%%%%%

      %%%%%%%%%%%%%%%%%%%%%%%%%%%%%%%%%%%%%%%%%%%%%%%%%%%%%%
      %%  TRUC AVEC MACHIN AU-DESSUS ET BIDULE AU-DESSOUS %%
      %%                (MODE MATHEMATIQUE)               %%
      %%%%%%%%%%%%%%%%%%%%%%%%%%%%%%%%%%%%%%%%%%%%%%%%%%%%%%

\def\build #1_#2^#3{\mathrel{\mathop{\kern 0pt #1}
     \limits_{#2}^{#3}}}

%%%%%%%%%%%%%%%%%%%%%%%%%%%%%%%%%%%%%%%%%%%%%%%%%%%%%%%%%%%%%%%%%%%

%%%%%%%%%%%%%%%%%%%%%%%%%%%%%%%%%%%%%%%%%%%%%%%%%%%%%%%%%%%%%%%%%%%%

         %%%%%%%%%%%%%%%%%%%%%%%%%%%%%%%%%%%%%%%%%%%%%%%%%%%
         %% FLECHES HORIZONTALES ET VERTICALES AVEC TRUC  %%
         %%   AU-DESSUS ET MACHIN AU-DESSOUS (MODE MATH)  %%
         %%%%%%%%%%%%%%%%%%%%%%%%%%%%%%%%%%%%%%%%%%%%%%%%%%%

\def\hfl[#1][#2][#3]#4#5{\kern -#1
 \sumdim=#2 \advance\sumdim by #1 \advance\sumdim by #3
  \smash{\mathop{\hbox to \sumdim {\rightarrowfill}}
   \limits^{\scriptstyle#4}_{\scriptstyle#5}}
    \kern-#3}

\def\vfl[#1][#2][#3]#4#5%
 {\sumdim=#1 \advance\sumdim by #2 \advance\sumdim by #3
  \setbox1=\hbox{$\left\downarrow\vbox to .5\sumdim {}\right.$}
   \setbox1=\hbox{\llap{$\scriptstyle #4$}\box1
    \rlap{$\scriptstyle #5$}}
     \vcenter{\kern -#1\box1\kern -#3}}

%%%%%%%%%%%%%%%%%%%%%%%%%%%%%%%%%%%%%%%%%%%%%%%%%%%%%%%%%%%%%%%%%%%%

                   %%%%%%%%%%%%%%%%%%%%%%%%%%%%%
                   %% DIAGRAMMES COMMUTATIFS  %%
                   %%%%%%%%%%%%%%%%%%%%%%%%%%%%%
\newdimen\sumdim
\def\diagram#1\enddiagram
    {\vcenter{\offinterlineskip
      \def\tvi{\vrule height 10pt depth 10pt width 0pt}
       \halign{&\tvi\kern 5pt\hfil$\displaystyle##$\hfil\kern 5pt
        \crcr #1\crcr}}}

%% SYNTAXE

%%  $$
%%  \diagram
%%   ... & & ... & \cr
%%   ... & & ... & \cr
%%   .................
%%   ... & & ... & \cr
%%  \enddiagram
%%  $$        

%%%%%%%%%%%%%%%%%%%%%%%%%%%%%%%%%%%%%%%%%%%%%%%%%%%%%%%%%%%%%%%%%%%%%%%%%%%%%%%%%%%%%%
%%%%%%%%%%%%%%%%%%%%%%%%%%%%%%%%%%%%%%%%%%%%%%%%%%%%%%%%%%%%%%%%%%%%%%%%%%%%%%%%%%%%%%

%\def\Doteq{\build{=}_{\scriptscriptstyle\bullet}^{\scriptscriptstyle\bullet}}

%%%%%%%%%%%%%%%%%%%%%%%%%%%%%%%%%%%%%%%%%%%%%%%%%%%%%%%%%%%%%%%%%%%%

               %%%%%%%%%%%%%%%%%%%%%%%%%%%%%%%%%%%
               %% ESPACE VERTICAL EN MODE MATH  %%
               %%%%%%%%%%%%%%%%%%%%%%%%%%%%%%%%%%%

\def\vspace[#1]{\noalign{\vskip #1}}

%%%%%%%%%%%%%%%%%%%%%%%%%%%%%%%%%%%%%%%%%%%%%%%%%%%%%%%%%%%%%%%%%%%%

                      %%%%%%%%%%%%%%%%%%%
                      %%  ENCADREMENT  %%
                      %%%%%%%%%%%%%%%%%%%

\def\boxit [#1]#2{\vbox{\hrule\hbox{\vrule
     \vbox spread #1{\vss\hbox spread #1{\hss #2\hss}\vss}%
      \vrule}\hrule}}

 %% INDIQUE LA FIN D'UNE PREUVE

%%%%%%%%%%%%%%%%%%%%%%%%%%%%%%%%%%%%%%%%%%%%%%%%%%%%%%%%%%%%%%%%%%%%

\catcode`\@=11
\def\Eqalign #1{\null\,\vcenter{\openup\jot\m@th\ialign{
\strut\hfil$\displaystyle{##}$&$\displaystyle{{}##}$\hfil
&&\quad\strut\hfil$\displaystyle{##}$&$\displaystyle{{}##}$
\hfil\crcr #1\crcr}}\,}
\catcode`\@=12

%%%%%%%%%%%%%%%%%%%%%%%%%%%%%%%%%%%%%%%%%%%%%%%%%%%%%%%%%%%%%%%%%%%%%%%%%%%%%%%%%%%%%%%%%%%%%%%%
%%%%%%%%%%%%%%%%%%%%%%%%%%%%%%%%%%%%%%%%%%%%%%%%%%%%%%%%%%%%%%%%%%%%%%%%%%%%%%%%%%%%%%%%%%%%%%

%\font\tengoth=eufm10            \font\tenbboard=msbm10
%\font\eightrm=cmr8              \font\eighti=cmmi8 
%\font\eightsy=cmsy8             \font\eightbf=cmbx8 
% \font\eighttt=cmtt8             \font\eightit=cmti8
% \font\eightsl=cmsl8             \font\eightgoth=eufm7
% \font\eightbboard=msbm7         \font\eightex=cmex10 at 8pt

%\font\eightgoth=eufm7           \font\sevenbboard=msbm7

% \font\sixrm=cmr6                \font\sixi=cmmi6
%\font\sixsy=cmsy6               \font\sixbf=cmbx6
% \font\sixgoth=eufm5 at 6pt      \font\sixbboard=msbm7 at 6pt

%\font\fivegoth=eufm5            \font\fivebboard=msbm5

\vsize= 17.0cm
\hsize= 13.0cm
\hoffset=3mm
\voffset=3mm

%%%%%%%%%%%%%%%%%
%% DEFINITIONS %%
%%%%%%%%%%%%%%%%%

\def\qed{\quad\raise -2pt\hbox{\vrule\vbox to 10pt{\hrule width 4pt 
\vfill\hrule}\vrule}} 

\def\qeda{\quad\raise -2pt\hbox{\vrule\vbox to 10pt{\hrule 
width 4pt height9pt 
\vfill\hrule}\vrule}}

\def\tend_#1^#2{\mathrel{
	\mathop{\kern 0pt\hbox to 1cm{\rightarrowfill}}
	\limits_{#1}^{#2}}}

\def\lfq{\leavevmode\raise.3ex\hbox{$\scriptscriptstyle
\langle\!\langle$}\thinspace}
\def\rfq{\leavevmode\thinspace\raise.3ex\hbox{$\scriptscriptstyle
\rangle\!\rangle$}}

\font\msamten=msam10
\font\msamseven=msam7
\newfam\msamfam

\textfont\msamfam=\msamten
\scriptfont\msamfam=\msamseven

\def\hexnumber #1%
{\ifcase #1 0\or 1\or 2\or 3\or 4\or 5
   \or 6\or 7\or 8\or 9\or A\or B\or C\or D\or E\or F
    \fi}

\mathchardef\leadsto="3\hexnumber\msamfam 20

\def\frac #1#2{{#1\over #2}}

%%%%%%%%%%%%%%%%%%%%%%%%
%% END OF DEFINITIONS %%
%%%%%%%%%%%%%%%%%%%%%%%%

%%%%%%%%%%%%%%%%%%%%%%%%%%%%%%%%%%%%%%%%% suites exactes %%%%%%%%%%%%%%%%%%%%%%%%%%%%%%%%%%%%%%%%%%%%%%%%%%%%%%%%%%%%%%

\def\rarr#1#2{\smash{\mathop{\hbox to .5in{\rightarrowfill}}
  \limits^{\scriptstyle#1}_{\scriptstyle#2}}} 

\def\larr#1#2{\smash{\mathop{\hbox to .5in{\leftarrowfill}}
  \limits^{\scriptstyle#1}_{\scriptstyle#2}}} 

\def\darr#1#2{\llap{$\scriptstyle #1$}\left\downarrow
  \vcenter to .5in{}\right.\rlap{$\scriptstyle #2$}}

\def\diagram#1{{\normallineskip=8pt
  \normalbaselineskip=0pt \matrix{#1}}}

%%%%%%%%%%%%%%%%%%%%%%%%%%%%%%%%%%%%%%%%%%%%%%%%%%%%%%%%%%%%%%%%%%%
%%%%%%%%%%%%%%%%%%%%%%%%%%%%%%%%%%%%%%%%%%%%%%%%%%%%%%%%%%%%%%%%%%%

%%%%%%%%%%%%%%%%%%%%%%%%%%%%
%%  AUTHOR ON EVEN PAGES  %%
%%%%%%%%%%%%%%%%%%%%%%%%%%%%

\auteurcourant{\eightbf Bodo Lass \hfill}

%%%%%%%%%%%%%%%%%%%%%%%%%%
%%  TITLE ON ODD PAGES  %%
%%%%%%%%%%%%%%%%%%%%%%%%%%

\titrecourant{\hfill\eightbf Une conjecture de Kontsevich et Shoikhet}

\vskip 1.5cm

%%%%%%%%%%%%%%%%%%%%%%%%%%%%%%%%%%%%%%%%%%%%%%%%%%%%%%%%%%%%%%%%%%%

%%%%%%%%%%%%%
%%  TITEL  %%
%%%%%%%%%%%%%

\centerline
{\fourteenbf
Une conjecture de Kontsevich et Shoikhet
}

\medskip

\centerline
{\fourteenbf
et la caract\'eristique d'Euler 
}

\vskip 15pt

%%%%%%%%%%%
%%  NAME %%
%%%%%%%%%%%  

\centerline{\bf BODO LASS}

\bigskip
\medskip

\centerline{\it A J.-P.~Jouanolou et A.~Al~Amrani,}

\smallskip

\centerline{\sevenrm dont j'ai oubli\'e l'invitation \`a cause de cette Note}

\bigskip
\medskip

%%%%%%%%%%%%%%%%
%%  ABSTRACT  %%
%%%%%%%%%%%%%%%%

{\parindent=.3in \narrower

\noindent
{\eightrm Salikhov~[8] a d\'emontr\'e une conjecture de Kontsevich et Shoikhet~[4] en la r\'eduisant \`a la consid\'eration de trois 
familles de graphes, une consid\'eration qui ne fut explicit\'ee que pour une de ces familles.
Nous montrons que cette conjecture de Kontsevich et Shoikhet n'est, en fait, qu'un cas tr\`es particulier du th\'eor\`eme classique 
sur la caract\'eristique d'Euler, bien explicit\'e par Teleman~[10].
}

\bigskip

\noindent
{\eightrm Salikhov~[8] has proved a conjecture of Kontsevich and Shoikhet~[4] by reducing it to the consideration of three 
families of graphs, a consideration which was left to the reader for two of those families.
We show, that the conjecture is just a very particular case of the classical theorem on the Euler characteristic,
well explicated by Teleman~[10].
}

\bigskip

\noindent 
{\eightit 1991 AMS Subject Classification:} {\eightrm 05C10, 05C25, 55U10}

\noindent
{\eightit Key words:} {\eightrm graph-complex, orientation, automorphism, Euler characteristic}

\par}

\vskip1.5pt

\parindent=.3in
\vskip 5mm

%%%%%%%%%%%%
%%  TEXT  %%
%%%%%%%%%%%%

\medskip

Les deux d\'efinitions et le th\'eor\`eme suivant constituent l'essentiel d'un article r\'ecent de Salikhov~:

\medskip

{\sc D\'efinition 1} ([3]). Soit $G = (V,E)$ un multigraphe, c'est-\`a-dire on accepte la pr\'esence de boucles et d'ar\^etes
multiples. Un automorphisme $P \in \hbox{\tenrm Aut}(G)$ induit une permutation~$|P_E|$ des ar\^etes de~$G$ et une transformation 
lin\'eaire non-d\'eg\'en\'er\'ee $P_E \in \hbox{\tenrm GL}(H_1(G,{\bb R}))$ sur le premier groupe d'homologie du graphe~$G$ 
(en tant qu'espace topologique) avec des coefficients r\'eels. 
L'homo\-morphisme d'orienta\-tion de Kontsevich $\Theta_K: \hbox{\tenrm Aut}(G) \to \{1,-1\}$ est d\'efini par la formule 
$\Theta_K(P) := \hbox{\tenrm sign}(|P_E|) \cdot \hbox{\tenrm sign}(\det[P_E,H_1(G,{\bb R})])$.

\medskip

{\sc D\'efinition 2} ([4]). Soit $G = (V,E)$ un multigraphe. Orientons (\`a l'aide de fl\`eches) les ar\^etes de~$G$ de fa\c con
arbitraire. Un automorphisme $P \in \hbox{\tenrm Aut}(G)$ induit une permutation~$P_V$ des sommets de~$G$ et une fonction 
$\epsilon_P: E \to \{1,-1\}$ d\'efinie par $\epsilon_P(e) := 1$ (resp.~$\epsilon_P(e) := -1$), si la permutation induite
$|P_E|$ respecte (resp.~renverse) l'orientation sur l'ar\^ete $e \in E$.
L'homo\-morphisme d'orientation de Shoikhet $\Theta_S: \hbox{\tenrm Aut}(G) \to \{1,-1\}$ est d\'efini par la formule
$\Theta_S(P) := \hbox{\tenrm sign}(P_V) \cdot \prod_{e \in E} \epsilon_P(e)$.

\medskip

Il est facile de montrer que le dernier produit est ind\'ependant de l'orientation choisie pour chaque ar\^ete, et l'on a

\medskip

{\sc Th\'eor\`eme} ([8]). {\it Soit $G = (V,E)$ un multigraphe connexe. Alors pour tout $P \in \hbox{\tenrm Aut}(G)$ on a }
$\Theta_K(P) = \Theta_S(P)$.

\medskip

Salikhov~[8] d\'emontre le th\'eor\`eme pr\'ec\'edent sur plus de deux pages en le r\'edui\-sant \`a la consid\'eration de trois 
familles de graphes, une consid\'eration qui ne fut explicit\'ee que pour une de ces familles. Selon~[8], ce th\'eor\`eme a \'et\'e 
conjectur\'e par Kontsevich et Shoikhet dans~[4], un article qui est situ\'e dans le contexte de~[1] et~[3]. De plus, Salikhov~[8] 
remercie {\lfq S.~V.~Duzhin for suggesting the problem and B.~Shoikhet for useful discussions\rfq}.
Nous tenons \`a souligner que nous n'avons pas eu la chance de voir l'article de Kontsevich et Shoikhet~[4] personnellement.

\bigskip

Avant de formuler les d\'efinitions ci-dessus, Salikhov~[8] a introduit {\lfq the so-called {\it half-edge language}\rfq} pour 
d\'ecrire les multigraphes, bien que les deux d\'efinitions fassent plut\^ot appel aux concepts classiques de la
th\'eorie des graphes et de la topologie alg\'ebrique (de Hopf ou Lefschetz). 

En fait, pour d\'emontrer la conjecture de Kontsevich et Shoikhet, il suffit de remplacer la trace par le d\'eterminant dans la 
d\'emonstra\-tion du th\'eor\`eme du point fixe de Lefschetz (dans le cas particulier des complexes un-dimensionnels, voir~[7], 
paragraphe 22).
C'est ce que nous proposons d'expliciter dans cette Note.

\bigskip

Soit $C_1(G,{\bb R})$ (resp.~$C_0(G,{\bb R})$) l'espace vectoriel engendr\'e par les ar\^etes (resp. sommets) du graphe~$G$.
Le leitmotif de la topologie alg\'ebrique est d'associer \`a un automorphisme le diagramme commutatif
\vskip-10pt
$$\diagram{
0 & \rarr{}{} & C_1(G,{\bb R}) & \rarr{\partial}{} & C_0(G,{\bb R}) & \rarr{\epsilon}{} & {\bb R} & \rarr{}{} & 0 \cr
  &  & \darr{P_E}{} & & \darr{}{P_V} &  & \darr{}{id}  \cr
0 & \rarr{}{} & C_1(G,{\bb R}) & \rarr{\partial}{} & C_0(G,{\bb R}) & \rarr{\epsilon}{} & {\bb R} & \rarr{}{} & 0 \cr
}$$
\vskip-5pt
\noindent
o\`u $P_E$ d\'esigne la permutation {\lfq sign\'ee\rfq} des ar\^etes orient\'ees selon la d\'efinition~2. Dans cette situation,
le polyn\^ome caract\'eristique est une application d'Euler-Poincar\'e d\'efinie sur le groupe d'Euler-Grothendieck (voir~[2] 
ainsi que~[5], chapitre~IV.3 et~XV.4, ou bien la nouvelle \'edition~[6], chapitre~XIV.3 et~XX.3). 

\eject

\noindent
En particulier, on a la relation suivante pour les d\'etermi\-nants des transformations induites sur les espaces d'homologie 
r\'eduits (voir~[10] ou~[9], paragraphe~64, exercice~28)~:
$$\openup2pt\eqalign{
        &\det[P_E,\tilde H_1(G,{\bb R})] \cdot \det[id,\tilde H_{-1}(G,{\bb R})] / \det[P_V,\tilde H_0(G,{\bb R})] \cr
\; = \; &\det[P_E,\,     C_1(G,{\bb R})] \cdot \det[id,\,     C_{-1}(G,{\bb R})] / \det[P_V,\,     C_0(G,{\bb R})].
}$$
\'Evidemment, $C_{-1}(G,{\bb R}) = {\bb R}$, $\tilde H_{-1}(G,{\bb R}) = 0$ et $\tilde H_0(G,{\bb R}) = 0$ puisque $G$ est connexe.
Par cons\'equent, 
$\det[id,C_{-1}(G,{\bb R})] \, = \, \det[id,\tilde H_{-1}(G,{\bb R})] \, = \, \det[P_V,\tilde H_0(G,{\bb R})] \, = \, 1$, 
d'o\`u
$$
\det[P_E,H_1(G,{\bb R})] \; = \; \det[P_E,\tilde H_1(G,{\bb R})] \; = \; 
\det[P_E,C_1(G,{\bb R})] / \det[P_V,C_0(G,{\bb R})].
$$
Mais $\det[P_E,C_1(G,{\bb R})] \, = \, \hbox{\tenrm sign}(|P_E|) \cdot \prod_{e \in E} \epsilon_P(e)$ et 
$\det[P_V,C_0(G,{\bb R})] \, = \, \hbox{\tenrm sign}(P_V)$. Ceci d\'emontre bien la conjecture de Kontsevich et Shoikhet.

\bigskip

{\it Remarque.}
On peut travailler partout avec des coefficients entiers (voir~[7], paragraphe 22). En particulier, on sait {\it a priori} 
que $\det[P_E,H_1(G,{\bb R})] \, \in {\bb Z}$, d'o\`u $\det[P_E,H_1(G,{\bb R})] \in \{1,-1\}$, puisque
$\det[P_E,H_1(G,{\bb R})] \cdot \det[P_E^{-1},H_1(G,{\bb R})] \, = \, 1$.
En outre, on pourrait utiliser l'homologie non-r\'eduite et supprimer l'hypoth\`ese {\lfq con\-nexe\rfq}
si l'on fait intervenir la signature de la permutation induite sur les composantes connexes de~$G$. 
Finalement, toutes sortes de g\'en\'eralisations de la conjecture de Kontsevich et Shoikhet aux endomorphismes (g\'en\'eralis\'es),
aux complexes d'une dimension sup\'erieure \`a un, etc. sont possibles. L'essentiel est l'applicabilit\'e du r\'esultat de
Teleman~[10].

\bigskip
\bigskip

{\it Remerciements.} Je voudrais remercier vivement D.~Foata et J.-P.~Jouanolou pour toute leur aide et assistance.

\vfill
\eject

\vglue 2cm

%%%%%%%%%%%%%%%%%%%%%%%
%%   BIBLIOGRAPHIE   %%
%%%%%%%%%%%%%%%%%%%%%%%

\centerline{\bf R\'ef\'erences bibliographiques}

\medskip

{\baselineskip=9pt\eightrm
\item{[1]} D.~B.~Fuchs, {\eightit Kogomologii beskonechnomernykh algebr Lie.} Nauka, 1984.
\smallskip 

\item{[2]} J.~L.~Kelley et E.~H.~Spanier, {\eightit Euler characteristics.} Pacific J.~Math.~{\eightbf 26} (1968), 317-339.
\smallskip 

\item{[3]} M.~Kontsevich, {\eightit Formal (non)-commutative symplectic geometry.}
           in: The Gelfand Mathematical Seminars 1990-1992, Birkh\"auser, 1993, 173-187.
\smallskip 

\item{[4]} M.~Kontsevich et B.~Shoikhet, {\eightit Formality conjecture, geometry of complex manifolds and 
           combinatorics of graph-complex.} preprint.
\smallskip 

\item{[5]} S.~Lang, {\eightit Algebra.} Addison-Wesley, 1965.
\smallskip 

\item{[6]} S.~Lang, {\eightit Algebra, Third Edition.} Addison-Wesley, 1997.
\smallskip 

\item{[7]} J.~R.~Munkres, {\eightit Elements of algebraic topology.} Addison-Wesley, 1984. 
\smallskip 

\item{[8]} K.~Salikhov, {\eightit On the orientation of graphs.} arXiv:math.CO/0012252, 5 pages.
\smallskip 

\item{[9]} G.~Scheja et U.~Storch, {\eightit Lehrbuch der Algebra, Teil 2.} B.~G.~Teubner, Stuttgart, 1988.
\smallskip 

\item{[10]} S.~Teleman, {\eightit Sur la formule d'Euler-Poincar\'e-Hopf.} Rev.~Math.~Pures Appl.~{\eightbf 2} (1957), 551-554.
\smallskip 
}

\bigskip

\vskip 12pt

{\eightrm

\rightline{\eightit Universit\'e Louis-Pasteur}
\rightline{\eightit D\'epartement de math\'ematique}
\rightline{\eightit 7, rue Ren\'e-Descartes}
\rightline{\eightit F-67084 Strasbourg}
\rightline{Courriel: lass@math.u-strasbg.fr} 

\vskip 6pt

\rightline{et}

\vskip 6pt

\rightline{\eightit RWTH Aachen}
\rightline{\eightit Lehrstuhl II f\"ur Mathematik}
\rightline{\eightit Templergraben 55} 
\rightline{\eightit D-52062 Aachen} 
\rightline{Courriel: lass@math2.rwth-aachen.de} 

}

\end